\documentclass[a4paper,11pt]{amsart}
\usepackage[backref=page, colorlinks=true, allcolors=blue]{hyperref}
\usepackage[english]{babel}

\usepackage{amsmath,amscd,amssymb,amsthm,amsxtra,enumerate}
\usepackage[foot]{amsaddr}
\usepackage{mathrsfs}
\usepackage{mathtools}
\usepackage{tikz}
\usepackage{color,graphicx}
\usepackage{epsfig,epstopdf}
\usepackage[margin=1.25in]{geometry}
\usepackage{enumerate}
\usepackage{hyperref}
\usepackage{subcaption}
\usepackage{wrapfig}
\usepackage{subcaption}
\usetikzlibrary{patterns}
\usetikzlibrary{decorations.pathreplacing}

\newcommand{\norm}[1]{\left\Vert#1\right\Vert}

\newcommand{\T}{\mathbb{T}}

\newcommand\reallywidehat[1]{\arraycolsep=0pt\relax%
\begin{array}{c}
\stretchto{
  \scaleto{
    \scalerel*[\widthof{\ensuremath{#1}}]{\kern-.5pt\bigwedge\kern-.5pt}
    {\rule[-\textheight/2]{1ex}{\textheight}} 
  }{\textheight} %
}{0.5ex}\\           
#1\\                 
\rule{-1ex}{0ex}
\end{array}
}

\theoremstyle{definition}

\numberwithin{equation}{section}

\allowdisplaybreaks

\author[A. Biswas]{Animesh Biswas}

\address{Department of Mathematics \\
{Missouri State University} \\
Springfield, MO, USA}
\email{ab7e@missouristate.edu}

\author[A. Huang]{Archie Huang}
\address{Department of Building, Civil and Environmental Engineering \\
Concordia University, Montreal, QC, Canada}
\email{archie.huang@concordia.ca}

\author[S. Agarwal]{Shaurya Agarwal}
\address{Department of Civil, Environmental \& Construction Engineering \\
{University of Central Florida} \\
Orlando, FL, USA}

\email{shaurya.agarwal@ucf.edu}

\author[C. Housholder]{Christopher Housholder}

\address{Department of Mathematics \\
{Missouri State University} \\
Springfield, MO, USA}
\email{clh284s@missouriState.edu}

\begin{document}

\title[Spatial-Temporal Nonlocal Traffic Dynamics]{Spatial-Temporal Nonlocal Traffic Dynamics: Analytical Properties, Adaptive Kernel Formulation, and Empirical Validation}
\maketitle

\begin{abstract}
This paper presents a new spatial‑temporal nonlocal traffic flow model formulated to overcome the boundedness limitations inherent in classical local formulations. The model introduces an adaptive kernel that captures both spatial and temporal nonlocal interactions, allowing the velocity at a given point to depend on aggregated downstream traffic conditions over a finite time horizon. This structure provides a more realistic representation of driver anticipation and reaction behavior. In addition to developing the model, we establish several key analytical properties that clarify the theoretical foundations of the proposed nonlocal framework. To assess its practical relevance, we conduct a detailed empirical validation using high‑resolution NGSIM trajectory data. The results demonstrate that the spatial‑temporal nonlocal model significantly improves the reconstruction of traffic density fields compared with traditional local macroscopic models, particularly in regimes where anticipation effects dominate. These findings highlight the potential of spatial‑temporal nonlocal traffic dynamics as a robust theoretical and data‑driven framework for capturing complex traffic behavior.
\end{abstract}


\section{INTRODUCTION}
Classical first‑order traffic models—most notably the LWR formulation \cite{Lighthill1955, whitham2011linear, Richards1956} $$\partial_t\rho(t,x)+\partial_x f(\rho(t,x))=0,$$ describe how the vehicle density $\rho: \mathbb{R}^+\times [a,b]\to[0,\rho_m]$ evolves along a road segment $(a,b)$ and where the values of the function output reside in $[0,\rho_m]$, $\rho_m$ being the jam density.  The flow function is defined by $f(\rho(t,x)) =\rho(t,x) v(\rho(t,x))$ linking density with the corresponding density‑dependent speed. 

A well‑known limitation of this model is that the velocity $v(\rho(t,x))$
 depends only on the density at the same point $\rho(t,x)$. Since real drivers react to downstream traffic, this purely local dependence fails to represent anticipatory behavior. To overcome that, 
 a nonlocal velocity model was used. Models of this type, such as those in \cite{blandin2016well} and \cite{goatin2016well}, replace the pointwise density in the velocity with a downstream convolution: $$\partial_t\rho(t,x)+\partial_x\left[ \rho(t,x)v(\rho(t,x)\ast_d \eta(x))\right]=0,$$ where $\ast_d$ is the \emph{downstream} convolution defined as $\rho(t,x)\ast_d \eta(x)=\int_x^{x+d}\rho(t,y)\eta(y-x)dy$.  The kernel $\eta(x)$ is a nondecreasing probability function supported on $[0,d]$ and it satisfies $\int^\infty_0 \eta(x)dx=1$.  This structure naturally incorporates driver look‑ahead by weighting downstream densities. The study in \cite{kachroo2023nonlocal} analyzed this nonlocal LWR model and validated its performance using NGSIM data. Because the model is nonlocal with a fixed kernel length, the authors employed modified boundary conditions, replacing standard boundary values with a thick boundary specification. In \cite{huang2024incorporating}, the framework was further extended to accommodate local boundary conditions through a variable‑length kernel, and the authors employed a physics‑informed neural network to compute solutions of the nonlocal LWR equation, illustrating the flexibility of nonlocal traffic modeling.

 Recent developments have significantly advanced the study of nonlocal formulations of traffic flow models and conservation laws \cite{du2012new, abreu2022lagrangian, treiber1999derivation, chiarello2021overview, bressan2020traffic, colombo2012class, piccoli2013transport, amorim2015numerical}. The motivation for these extensions is clear: the physical behavior of drivers inherently reflects an anticipation of downstream conditions rather than reliance solely on the immediate local environment. This anticipatory behavior corresponds to a spatial averaging of traffic density ahead of the vehicle, underscoring the importance of nonlocal modeling approaches in accurately representing realistic traffic dynamics. Beyond spatial nonlocality, a smaller body of work has also explored temporal nonlocal effects. In particular, \cite{du2022space} introduced a framework in which the density at a given location and time is influenced by the density at a downstream position at an earlier time. Their study established the well‑posedness of the resulting model and analyzed the associated conservation law for this spatial‑temporal nonlocal LWR formulation. 

 In the broader context, it is worth noting that beyond transportation engineering, nonlocal models have garnered significant attention within the theoretical and applied mathematics communities as well as various engineering disciplines. Their appeal stems from the ability of nonlocal formulations to represent multi‑scale interactions that cannot be captured adequately by purely local models. As a result, nonlocal frameworks have been successfully employed across a wide range of physical phenomena, including semi‑permeable membrane processes in cellular biology \cite{biswas2021harnack, biswas2021regularity}, nonlocal curvature and perimeter functionals in image processing \cite{biswas2026nonlocal}, and dynamic fracture mechanics \cite{silling2000reformulation}.


\vspace{1mm}
\noindent \textbf{Contributions:}  
This paper makes the following contributions: 
\begin{itemize} 
\item It provides a detailed investigation of the spatial–temporal dependence of the traffic density variable, examining how different parameter regimes influence the resulting dynamics. 
\item It provides the key analytical properties of the proposed spatial-temporal nonlocal model, including well‑posedness considerations and qualitative behavior under various kernel choices. 
\item It introduces an adapted Lax–Friedrichs numerical scheme tailored to the time-dependent nonlocal structure, ensuring stable and accurate approximate solutions. 
\item {It performs a rigorous empirical validation using high‑resolution NGSIM trajectory data, demonstrating that the proposed spatial‑temporal nonlocal model significantly improves the reconstruction of traffic density fields relative to traditional local macroscopic models.}
\item It analyzes the influence of kernel length and temporal dependence on boundary behavior, offering insights into how adaptive kernels affect boundary conditions and model calibration in practical implementations.
\item {It provides a comprehensive study of nonlocal boundary value problems arising in traffic applications, comparing multiple techniques for handling boundary effects in spatial‑temporal nonlocal flow models.}
\end{itemize}

\textbf{Outline:} Section \ref{sec:classical} reviews the classical (local) LWR model and the associated fundamental diagrams, including those of Greenshields and Underwood. Section \ref{sec:nonlocal} reviews the space-dependent nonlocal LWR formulation. Section \ref{sec:kernel} develops the analytical framework for the spatial-temporal nonlocal LWR model for both fixed and variable‑length kernels. Section \ref{sec:eval} assesses the effectiveness of the proposed nonlocal model and kernel choices using realistic traffic data, and Section \ref{sec:conc} concludes the paper with a summary of findings and directions for future developments in nonlocal traffic modeling.

\section{Classical LWR model}\label{sec:classical}
This section summarizes the classical (local) LWR traffic flow model and recalls two widely used fundamental diagrams: those of Greenshields and Underwood.
\subsection{The LWR Conservation Law}
In the Eulerian framework, with spatial–temporal coordinates $\mathbf{X}=(t, \, x)$, the Lighthill–Whitham–Richards (LWR) model is expressed as the conservation of vehicles:
\begin{equation} \label{eqn:lwr_conservation_recall}
{\mathrm{For} \,\, (t, x) \, \in \, \mathbb{R} \, \times \, \mathbb{R}^{+} \, :   \,\,\,  \frac{\partial \rho(t, \, x)}{\partial t} + \frac{\partial f(t, \, x)}{\partial x} = 0}
\end{equation}
Here $\rho(t,x)$ denoted the traffic density and $f(t,x)$ denoted the corresponding traffic flow.

\subsection{Greenshields Fundamental Diagram}
Fundamental diagrams specify functional relationships among traffic variables. The Greenshields model \cite{Greenshields1935} postulates a linear dependence between speed and density. With jam density $\rho_m$
 and free‑flow speed $v_f$, the flow–density and speed–density relations are:
\begin{equation} \label{eqn:greenshields_fundamental_diagram_recall}
\begin{cases}
     f(\rho(t, x)) = \rho(t,  x) \, v_f \bigg(1 - \frac{\rho(t, \, x)}{\rho_m}\bigg)
            \\[2pt]
            v(\rho(t, \, x))   =  v_f \left(1 - \frac{\rho(t, \, x)}{\rho_m} \right)
\end{cases}       
\end{equation}
Substituting the Greenshields relation into the conservation law yields a scalar PDE in terms of density alone:
\begin{equation} \label{eqn:lwr_conservation_density_recall}
  {\frac{\partial \rho(t, x)}{\partial t} + v_f \left(1 - \frac{2\rho(t, x)}{\rho_m}\right)\frac{\partial \rho(t, x)}{\partial x}   = 0}
\end{equation}


%
\subsection{Underwood Fundamental Diagram}

An alternative representation is provided by the Underwood model \cite{underwood1961speed}, which assumes an exponential decay of speed with density. With critical density $\rho_C$, the relations are:
\begin{equation} \label{eqn:underwood_fd}    
     \begin{cases}
         f(\rho(t, \, x))  =  \rho(t, \, x) \, v_f \, e^{- \frac{\rho(t, \, x)}{\rho_c}}
            \\[2pt]
            v(\rho(t, \, x))   =  v_f \, e^{- \frac{\rho(t, \, x)}{\rho_c}}
      
     \end{cases}     
\end{equation}

Combining the Underwood diagram with the LWR conservation law again produces a density‑based PDE:
\begin{equation} \label{eqn:lwr_density_underwood}
  {\frac{\partial \rho(t, x)}{\partial t} + v_f \, e^{- \frac{\rho(t, \, x)}{\rho_c}} \, \left(1 - \frac{\rho(t, x)}{\rho_c}\right)\frac{\partial \rho(t, x)}{\partial x}   = 0}
\end{equation}

\section{Nonlocal space-dependent traffic model}\label{sec:nonlocal}
\begin{figure}[htbp!]
\centering
\centering
\begin{tikzpicture}
\draw[line width=1mm] (-1,1)--(6.25,1);
\draw[dashed] (-1,0)--(6.25,0);
\draw[line width=1mm] (-1,-1)--(6.25,-1);
\node at (0,0) (c) {};
\shade[shading=radial, inner color=blue!100, outer color=blue!10] (c) -- (-25:2) arc(-25:25:2) -- cycle;
\draw[very thick] (c) -- (-25:2) arc(-25:25:2) -- cycle;
\draw[-latex,line width=0.75mm] (c.center)--+(1.5,0) node[below] {$v(\rho_d)$};
\node at (0,0) (c) {};
{\includegraphics[scale=0.05,angle=0]{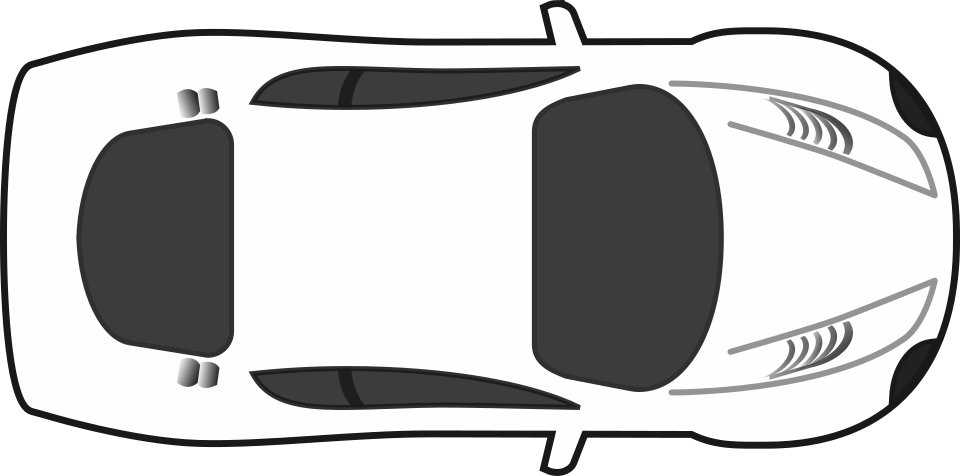}};

\node at (-1,-1.5) {$x=a$};
\draw[|-|] (0.95,0.75)--(2.2,0.75);
\node at (0.60,0.75) {$\eta(x)$};
\node at (0.75,-1.5) {$x$};
\node at (2,-1.5) {$x+d$};
\node at (6.25,-1.5) {$x=b$};
\node at (4.25,-0.5) {$\rho_d(t,x)$};
\end{tikzpicture}
\caption{Nonlocal Traffic Density Field}\label{Fig:Back:LWRNLtr}
\end{figure}
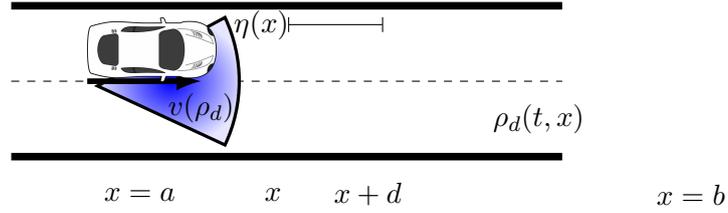

The nonlocal model considered here (Figure~\ref{Fig:Back:LWRNLtr}) describes traffic density $\rho(t,x) \in [0, \rho_m]$ on the interval $x \in (0, d)$. The look‑ahead density $\rho_d(t,x)$ can be obtained by convolving $\rho(t,x)$ with a non-increasing kernel function $\eta(x)$, that generalizes the probability density function $\eta(x)$, defined earlier.

The variable $\rho_d(t,x)$ is the smoothed look-ahead traffic density obtained by a convolution of $\rho(t,x)$ with a look-ahead kernel function $\eta(x)$. The kernel function $\eta(x)$ generalizes the probability density function $\eta(x)$, defined earlier. 

\begin{equation}\label{eq:Fig:Back:nlLWRRhod}
\begin{gathered}
\rho_d(t,x)=\rho(t,x) {\ast} \eta(x)
=\int_x^{x+d}\rho(t,y)\eta(y-x)dy
\end{gathered}
\end{equation}
The kernel satisfies the normalization condition
\begin{equation}\label{eq:Fig:Back:nlLWReta}
\int_0^d \eta(x)\,dx = 1 .
\end{equation}

The resulting nonlocal conservation law is
\begin{equation}\label{eq:Fig:Back:nlLWR}
\partial_t\rho(t,x)+\partial_x\!\left[\rho(t,x)\,v(\rho_d)\right]=0 .
\end{equation}
For the velocity function $v(\rho_d)$, any decreasing relation may be used; here we adopt the Greenshields‑type form
\begin{equation}\label{eq:Fig:Back:nlLWRGr}
v(\rho_d)=v_f\!\left(1-\frac{\rho_d}{\rho_m}\right),
\end{equation}
where $v_f$ is the free‑flow speed and $\rho_m$ the jam density.

\subsection{A General Initial Boundary Value Problem}\label{sec:Back:GM}
For completeness, we briefly summarize the weak formulation and well‑posedness framework for the space‑dependent nonlocal conservation law. Full details can be found in \cite{de2017initial} and the traffic‑flow adaptation in \cite{kachroo2023nonlocal}. The model is posed as the following initial–boundary value problem:
\begin{equation}\label{eq:Fig:Back:nl-GM}
\begin{gathered}
\partial_t\rho(t,x)+\partial_xf(t,x,\rho,\rho(t,x)*\eta(x))=0,\quad (t,x)\in\mathbb{R}^+\times(a,b)\\
\rho(0,x)=\rho_0,\quad x\in(a,b)\\
\rho(t,a)=\rho_a(t),\quad t\in\mathbb{R}^+\\
\rho(t,b)=\rho_b(t),\quad t\in\mathbb{R}^+\\
\end{gathered}
\end{equation}
where $f\in\mathbf{C}^2(\mathbb{R}^+\times[a,b]\times\mathbb{R}\times\mathbb{R};\mathbb{R})$ satisfies
\begin{equation}\label{eq:Fig:Back:nl-GMf}
\begin{gathered}
f(t,x,0,r)=0,\quad \forall t,x,r,\\
\sup_{t,x,\rho,r}|{\partial_{\rho}}f(t,x,\rho,r)|<K_1,\\
\sup_{t,x,r}|\partial_xf(t,x,\rho,r)|<K_2|\rho|, \\
\sup_{t,x,r}|\partial_rf(t,x,\rho,r)|<K_2|\rho|,\\
\sup_{t,x,r}|\partial_{xx}f(t,x,\rho,r)|<K_2|\rho|, \\
\sup_{t,x,r}|\partial_{xr}f(t,x,\rho,r)|<K_2|\rho|, \\
\sup_{t,x,r}|\partial_{rr}f(t,x,\rho,r)|<K_2|\rho|\\
\end{gathered}
\end{equation}
We have $K_1>0$, $K_2>0$ which are constants, and $\eta\in(\mathbf{C}^1\cap\mathbf{W}^{1,\infty})(\mathbb{R}, \mathbb{R})$ with
\begin{equation}\label{eq:Fig:Back:eta}
\begin{gathered}
\int_{\mathbb{R}}\eta(x)dx=1
\end{gathered}
\end{equation}

Because solutions may develop discontinuities, the problem is interpreted in the Kruzhkov entropy sense. Entropy inequalities incorporate both initial and boundary data, and boundary conditions are enforced through entropy‑admissible traces. Under the stated assumptions on $f, \eta$, and the data, the general theory guarantees a unique entropy solution $\rho \in BV$. 
 The theorem further establishes stability estimates: uniform bounds on the $L^\infty, L^1$, and total variation norms of the solution, continuity in time in the $L^1$ sense, and continuous dependence on initial and boundary data. These properties collectively ensure existence, uniqueness, and stability of solutions for the nonlocal LWR model.

\subsection{Kernel function}
We now revisit several kernel formulations previously developed for computing the space‑dependent nonlocal traffic density within the LWR framework. A comprehensive treatment of these kernels is provided in \cite{huang2024incorporating}, where the authors also integrate the nonlocal LWR model into a PIDL framework, demonstrating how nonlocal traffic dynamics can be combined with data‑driven learning approach.
Here we give a brief excerpt of that analysis. Two fixed‑length kernels are first presented: a constant kernel that averages downstream density uniformly, and a linearly decreasing kernel that assigns greater weight to nearby traffic conditions. 
{
\begin{equation} \label{eqn:constant_kernel}
\eta(x) =
\begin{cases}
    \frac{1}{d}, \;\;\; x \in [0, \;  d] \\
    0, \quad \text{else}    
\end{cases}
\end{equation}
}

\begin{figure}[htbp]
    \hspace{-5mm}\centerline{\includegraphics[width=0.3\textwidth]{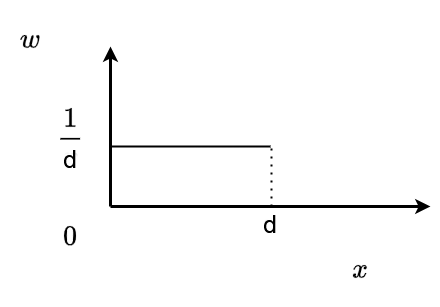}}
    \caption{The Constant Kernel}
    \label{fig:constant_kernel}
\end{figure}

{
\begin{equation} \label{eqn:linearly_decreasing_kernel}
\eta(x) =
\begin{cases}
     \frac{2}{d}(1 - \frac{x}{d}), \;\;\; x \in [0, \; d] \\
     0, \quad \text{else}  
\end{cases}    
\end{equation}
}

\begin{figure}[htbp]
    \hspace{-5mm}\centerline{\includegraphics[width=0.3\textwidth]{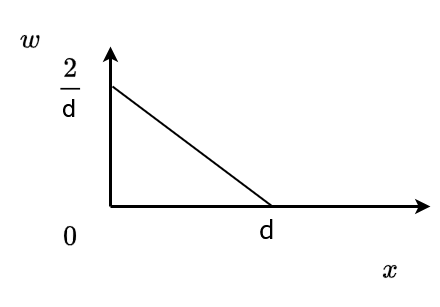}}
    \caption{The Linearly Decreasing Kernel 
    }
    \label{fig:linearly_decreasing_kernel}
\end{figure}

To address situations where thick downstream data are not available, \cite{huang2024incorporating} introduces variable‑length kernels, whose support automatically shrinks as a point approaches the boundary. In practice, this means the kernel adjusts to the available downstream domain, using a shorter look‑ahead window whenever the full kernel would extend beyond the boundary. This allows the nonlocal density to be computed even when only thin‑boundary data are available, for example, at a signalized intersection where measurements exist at the stop line but not downstream.

\section{Nonlocal time dependent traffic model} \label{sec:kernel}
We now turn to the central contribution of this work, namely the development of a spatial‑temporal nonlocal traffic model that captures both look‑ahead interactions and delayed driver responses within a unified LWR‑type framework.
It is important to note that, beyond spatial nonlocality, a smaller but growing body of work has also begun to investigate temporal nonlocal effects, which introduce an additional layer of realism by acknowledging that driver responses are not instantaneous. In real traffic, drivers react not only to the current downstream density but also to conditions they observed moments earlier, reflecting finite reaction times and the time required to process and act on visual information. Hence, it is useful to incorporate
time delays of this traffic density information in the distance. Incorporating such delays leads to models in which the flux at time $t$ depends on density values at earlier times $t-\tau$, or more generally on a weighted history of past states. These temporal memory effects fundamentally alter the structure of the conservation law, transforming it from a purely instantaneous relation into one with spatial‑temporal coupling.

In particular, \cite{du2022space} proposed a formulation in which the density at a given point and time is influenced by the density at a downstream location evaluated at a previous time. In particular, the nonlocal conservation law is given by the following equation:
\begin{align}\label{eq:nonlocal_space_time_conservation}
    \partial_t\rho(t,x)+\partial_x (\rho(t,x) v(\rho_d(t,x)))=0,
\end{align}
or,
\begin{align}\label{eq:nonlocal_space_time_conservation_flux}
    \partial_t\rho(t,x)+\partial_x (f(t,x, \rho, \rho_d))=0,
\end{align}
where
\begin{equation}\label{eq:nonlocal_density_time_space}
   \rho_d(t,x) = \int^\infty_0 \rho(t-\gamma s, x+s) \eta(s) ds 
\end{equation}
and $f(t,x, \rho, \rho_d) = \rho(t,x) v(\rho_d(t,x))$, $\eta(s)$ is a non-negative, non-increasing kernel function such that $\int^\infty_0 \eta(s) ds = 1$ and $\gamma$ is a propagation delay. If we assume that there is no propagation delay, i.e $\gamma = 0$, then we get back the spatially dependent nonlocal model, i.e \eqref{eq:Fig:Back:nlLWRRhod}. 

\begin{figure}
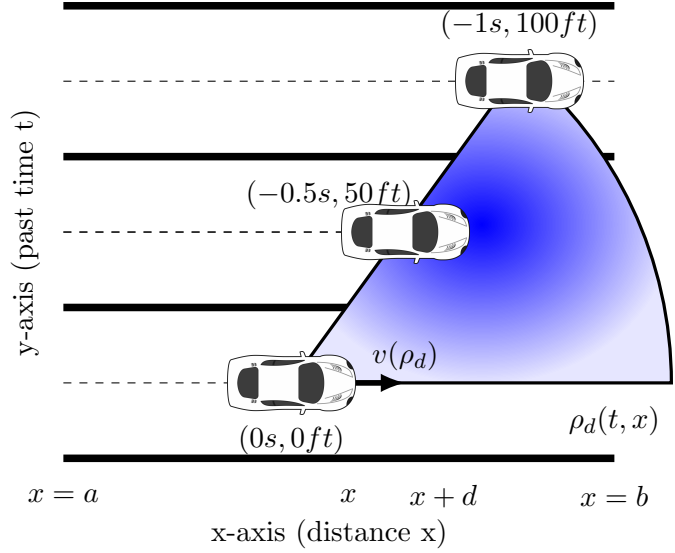

\begin{tikzpicture}

\newcommand{\Road}[5]{%
  \draw[line width=1mm] (#4,#1) -- (#5,#1);
  \draw[line width=1mm] (#4,#3) -- (#5,#3);
  \foreach \yy in {#2} { \draw[dashed] (#4,\yy) -- (#5,\yy); }%
}

\newcommand{\CarWedge}[1]{%
  \coordinate (carCenter) at (0,0);
  \shade[shading=radial, inner color=blue!100, outer color=blue!10]
    (carCenter) -- (0:5) arc(0:50:5.5) -- cycle;
  \draw[very thick]
    (carCenter) -- (0:5) arc(0:50:5.5) -- cycle;
  \draw[-latex,line width=0.75mm]
    (carCenter) -- +(1.5,0) node[above] {$v(\rho_d)$};
  \node at (carCenter) {\includegraphics[scale=#1,angle=0]{new.png}};
}

\def\xL{-1}
\def\xR{6.25}

\Road{1}{0}{-1}{\xL}{\xR}
\Road{3}{2,0,-2}{-3}{\xL}{\xR}

\node at (\xL,-3.5) {$x=a$};
\node at (\xR,-3.5) {$x=b$};

\node at (5,2.75) {$(-1s, 100ft)$};
\node at (2.5,0.5) {$(-0.5s, 50ft)$};
\node at (2,-2.75) {$(0s,0ft)$};

\node[rotate=90] at (-1.5,0) {y-axis (past time t)};
\node at (2.5,-4) {x-axis (distance x)};

\begin{scope}[shift={(2,-2)}]
  \draw[|-|] (0.75,0.75) -- (2,0.75);
  \node at (0.75,-1.5) {$x$};
  \node at (2,-1.5) {$x+d$};
  \node at (4.25,-0.5) {$\rho_d(t,x)$};
  \CarWedge{0.05}
\end{scope}

\node at (3.5,0) {\includegraphics[scale=0.05,angle=0]{new.png}};
\node at (5,2) {\includegraphics[scale=0.05,angle=0]{new.png}};

\end{tikzpicture}
\caption{Spatial-Temporal nonlocal traffic density field}
\end{figure}

This construction captures the idea that drivers anticipate future traffic conditions based on information perceived earlier, effectively blending spatial look‑ahead with temporal delay. Their analysis demonstrated that such a model remains mathematically well‑posed despite the introduction of delay terms, and they established existence, uniqueness, and stability of solutions for the resulting spatial‑temporal nonlocal LWR equation. Moreover, the presence of temporal nonlocality introduces new dynamical behaviors that differ significantly from those in classical or purely spatially nonlocal models. The work of \cite{du2022space} thus provides a foundational step toward understanding how memory and delayed anticipation shape macroscopic traffic flow.
We begin by reviewing the theoretical framework for the space–time nonlocal model presented in \cite{du2022space}. 
\subsection{Assumptions for the theoretical model}
In \cite{du2022space}, the authors made certain assumptions to prove the well-posedness of the initial value problems. 
\begin{itemize}
\item They assumed that the maximum density $\rho_m=1$.
    \item The velocity variable, as a function of the density is second order differentiable  and strictly decreasing in $[0, 1)$. We can see that the Greenshields model satisfies the assumption, if we assume use the normalized density to be $\frac{\rho}{\rho_m}$.
    \item The kernel function, $\eta(s)$, in addition to being non-negative, non-increasing and integrability condition, is $C^1([0,\infty))$ and satisfies,
    $$ \eta'(s) \leq -\beta \eta(s), \quad \forall s>0,$$
    where $\beta >0$ is a constant. 
\end{itemize}
They defined the $\rho_0(x)$ as initial data function and then defined
$$\rho_-(t,x) = \rho_0(x), \quad \text{for all}~(t,x) \in (-\infty, 0) \times \mathbb{R},$$
to be the past-time data which was used to compute the nonlocal density. Notice that, since the kernel function is defined over $[0,\infty)$, we need all past time data for this problem, which the authors created just by extending one set of initial data for all time $t<0$. Then they proved, given 
$$\gamma \leq \gamma_{max} = \min \bigg\{\frac{1}{3(v_f + \norm{v'}_{\infty})}, \frac{\beta}{\eta(0)\norm{v'}_{\infty}}  \bigg\},$$
where $\norm{v'}_{\infty}$ denotes the maximum value of $v'$, \cite{du2022space}[Theorem 1.2], that the following initial value problem has a weak solution to \eqref{eq:nonlocal_space_time_conservation} given the initial condition $\rho_{-}(t,x)$.

\subsection{Nonlocal boundary data}
Next we examine a key aspect of the nonlocal boundary‑value problem: namely, determining the appropriate form of boundary data for the space–time nonlocal model. 
In the space–time nonlocal model, the formulation naturally requires both thick boundary data and thick initial data. This is a key distinction from the purely spatial nonlocal model. In \cite{kachroo2023nonlocal}, the author considers only the spatially nonlocal setting and therefore imposes a thick boundary condition in space. For the space–time nonlocal problem, however, one must also prescribe thick initial data, which effectively acts as part of the thick boundary data in the space–time domain. As shown in the theoretical analysis of \cite{du2022space}, the authors define a function $\rho_{-}(t,x)$ that extends the initial profile $\rho_0(x)$ for all $t< 0$.

\begin{figure}[h!]
\begin{tikzpicture}[scale=2] 
\def\inner{1.0} 
\def \outer{1.3} 
\draw[very thick] plot [smooth cycle, tension=0.9] coordinates { (0, \outer) (\outer, 0.1) (0.7, -\outer) (-0.8, -0.9) (-\outer, 0.2) }; 
\begin{scope} 
 \fill[blue!15] plot [smooth cycle, tension=0.9] coordinates { (0, \outer) (\outer, 0.1) (0.7, -\outer) (-0.8, -0.9) (-\outer, 0.2) }; 
\fill[green!15] plot [smooth cycle, tension=0.9] coordinates { (0, \inner) (\inner, 0.1) (0.6, -\inner) (-0.6, -0.7) (-\inner, 0.2) }; 
\end{scope}
\draw[thick] plot [smooth cycle, tension=0.9] coordinates { (0, \outer) (\outer, 0.1) (0.7, -\outer) (-0.8, -0.9) (-\outer, 0.2) }; 
\draw[very thick] plot [smooth cycle, tension=0.9] coordinates { (0, \inner) (\inner, 0.1) (0.6, -\inner) (-0.6, -0.7) (-\inner, 0.2) }; 
\node at (0,0) {$\Omega$}; \node[blue!60!black] at (1.1,0.3) {$d$};
\draw[<->] (1.05,0.15) -- (1.2,0.25); 
\node at (0, \outer-0.12) {$\Gamma_d$};
\end{tikzpicture}
\caption{Nonlocal boundary data.}
\label{fig:nonlocal_bdd_general}
\end{figure}
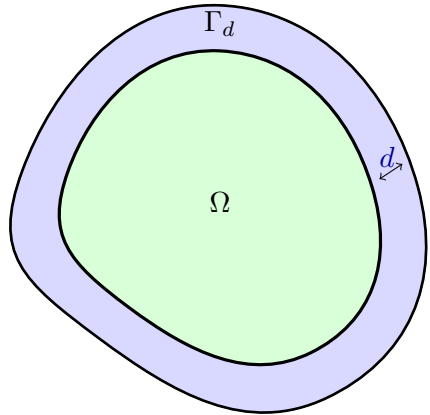

For our nonlocal traffic problem the domain is bounded. Specifically, we consider $(0,T) \times (a,b)$ as our two‑dimensional space–time domain. The nonlocal operator requires data not only inside the domain but also in a surrounding collar region in space–time.  In a prototypical nonlocal problem in a bounded domain, we consider $\Omega \subset \mathbb{R}^2$ to be the domain. We assume $\mu(x,y)$ to be the kernel function supported in the ball $B_d(0)$, centered on the point $(0,0)$ and radius $d$. To define a nonlocal quantity via convolution, one needs values not only in $\Omega$, rather in the enlarged set $\Omega_d$, which is defined as:
$$\Omega_d = \Omega \cup \{(x_0, y_0) \in  \mathbb{R}^2 |\quad  dist((x_0, y_0), U) < d \},$$
where the distance from any point $(x_0,y_0)$ to the open set $\Omega$ is defined as,
$$ dist((x_0, y_0), U) = \inf_{(x,y) \in \Omega} \{ \sqrt{(x-x_0)^2 + (y-y_0)^2} \}.$$
The associated collar region is then defined as, see Fig \ref{fig:nonlocal_bdd_general}, 
$$\Gamma_d = \{(x_0, y_0) \in  \mathbb{R}^2 |\quad  dist((x_0, y_0), U) < d \}.$$
It is also worth noting that for many nonlocal problems it is natural, and sometimes essential, to take $d=\infty$,in which case the collar region becomes unbounded. For further discussion and additional examples, the interested reader may consult the references \cite{biswas2021regularity, biswas2026nonlocal, BISWAS2024110474,   biswas2021harnack, foss2019bridging}.

In the space–time setting, the thick initial data plays the role of this collar in the temporal direction, ensuring that the nonlocal operator is well defined throughout the domain. In the traffic flow problem, it is important to emphasize that the nonlocal interaction in time is strictly retrospective: the kernel samples only past values and never future ones. Consequently, the temporal collar is required only near the initial time, not at the terminal time $t=T$. In contrast, the spatial interaction is look‑ahead in the $x$-direction  meaning that each point interacts with points located ahead of it (toward increasing $x$). Therefore, the spatial collar is needed only near the right boundary $x=b$, and not near the left boundary $x=a$. Also, we note that if the thick boundary in the space is of length $d$ then the thick boundary in the time is of length $\gamma d$.

We now turn to the next contribution of this work, where we propose two distinct formulations of thick boundary conditions for the space–time nonlocal model. In all earlier studies of the purely spatial nonlocal case \cite{kachroo2023nonlocal}, \cite{huang2024incorporating}, the thick‑boundary treatment has always followed an approach which is equivalent to the second approach described in the following. 
\subsubsection{Continuously extended data}\label{sub:cont_bdd}
\begin{figure}[h!]
\begin{tikzpicture}[scale=1.4]

\def\T{5} 
\def\gd{0.5} 
\def\a{0} 
\def\b{4} 
\def\d{1} 

\draw[thick] (0, \a) rectangle (\T, \b+\d); 
\fill[red!35] (0, \a) rectangle (\gd, \b+\d); 
\fill[blue!35] (\gd, \b+\d) rectangle (\T, \b); 
\fill[green!15] (\gd,\a) rectangle (\T, \b);

\fill[blue!15] (0, \b+\d) rectangle (\gd, \b);

\draw[very thick] (\gd,\a) rectangle (\T, \b);

\node at (2.0,2.5) {$\text{Interior Domain}$};
\node at (2.0,2) {$(0, T)\times(a,b)$};
\node[blue!60!black] at (2.0, 4.6) {Spatial collar of width $d$}; 
\node[red!60!black, rotate=90] at (0.15,2.0) {Temporal collar of width $\gamma d$}; 
\node at (-0.3,2) {$x$};
\node at (2.5,-0.3) {$t$}; 
\node at (-0.3,0) {$a$}; 
\node at (-0.3,\b+\d) {$b+d$}; 
\node at (-0.3,\b) {$b$}; 
\node at (\gd,-0.2) {$0$}; 
\node at (0,-0.3) {$-\gamma d$}; 
\node at (\T,-0.2) {$T$};
\end{tikzpicture}
\caption{Nonlocal `thick' boundary by extending local boundary data.}
\label{fig:Du_bdd}
\end{figure}
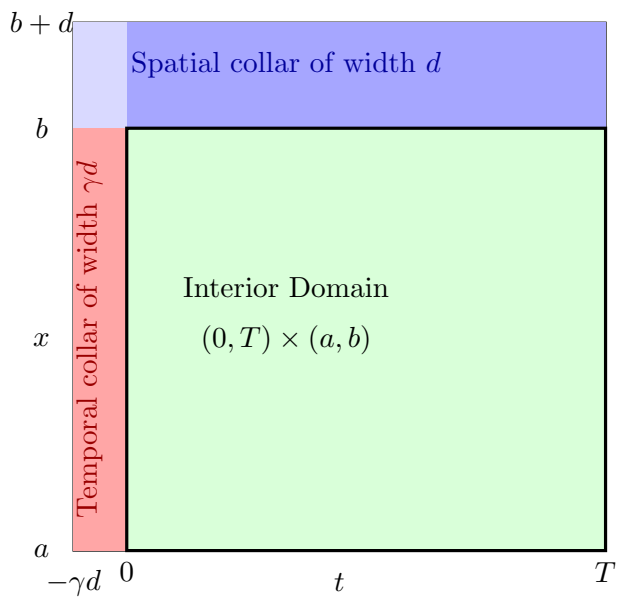

Our first formulation instead constructs the thick boundary region by continuously extending the boundary data backward in time near $t=0$ and forward in space near $x=b$ ensuring that the nonlocal operator remains well defined throughout the enlarged collar region. 
If we extend the boundary data in this manner, we obtain a construction that is conceptually similar to the extension method used in \cite{du2022space}, as illustrated in Figure~\ref{fig:Du_bdd}.

In the red region, i.e, in $[\gamma d, 0] \times [a,b]$, we use 
$$u(t,x) = u(0,x), \quad \text{for all}~ x \in [a,b], \quad t \in (-\gamma d, 0).$$
On the other hand, in the dark blue region, i.e in $[0, T] \times [b, b+d]$ we use,
$$u(t,x) = u(t,b), \quad \text{for all}~ x \in [b, b+d], \quad t \in (0, T).$$
In the light blue region, i.e in $[-\gamma d, 0] \times [b, b+d]$ we use 
$$u(t,x)=u(0,b), \quad \text{for}~(t,x) \in [-\gamma d, 0] \times [b, b+d]. $$
We notice that the last extension is well-defined and maintains continuity with the other two extensions.

\subsubsection{Known thick boundary}\label{sub:known_bdd}
The second formulation assumes that thick boundary information is already available on a prescribed subset of the space–time domain.
\begin{figure}[h!]
\begin{tikzpicture}[scale=1.4]

\def\T{5} 
\def\gd{0.5} 
\def\a{0} 
\def\b{4} 
\def\d{1} 


\draw[thick] (0, \a) rectangle (\T, \b); 
\fill[red!35] (0, \a) rectangle (\gd, \b-\d); 
\fill[blue!35] (\gd, \b- \d) rectangle (\T, \b); 
\fill[green!15] (\gd,\a) rectangle (\T, \b-\d);

\fill[blue!15] (0, \b-\d) rectangle (\gd, \b);

\draw[very thick] (\gd,\a) rectangle (\T, \b-\d);

\node at (2.0,2.5) {$\text{Interior Domain}$};
\node at (2.0,2) {$(\gamma d, T)\times(a,b-d)$};
\node[blue!60!black] at (2.0,3.4) {Spatial collar of width $d$}; 
\node[red!60!black, rotate=90] at (0.15,2.0) {Temporal collar of width $\gamma d$}; 
\node at (-0.3,2) {$x$};
\node at (2.5,-0.3) {$t$}; 
\node at (-0.3,0) {$a$}; 
\node at (-0.3,\b-\d) {$b-d$}; 
\node at (-0.3,\b) {$b$}; 
\node at (0,-0.2) {$0$}; 
\node at (\gd,-0.3) {$\gamma d$}; 
\node at (\T,-0.2) {$T$};
\end{tikzpicture}
\caption{Nonlocal `thick' boundary with known data }
\label{fig:thick_bdd}
\end{figure}
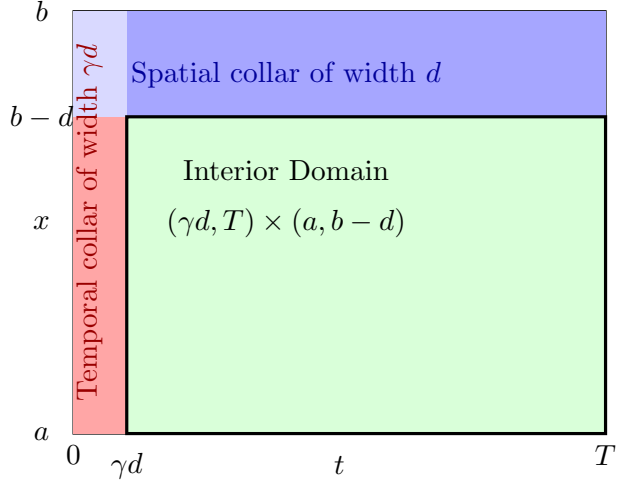

As one can see in Figure \ref{fig:thick_bdd}, the boundary data are indeed known in a `thick' region. In this formulation, the solution is known on the set $[0, \gamma d] \times [a, b] \cup [\gamma d, T] \times [b-d, d]$. 
As a result, the effective computational domain where the nonlocal operator can be evaluated without requiring data outside the prescribed region is the rectangle $(\gamma d, T) \times (a, b-d)$.

\section{Numerical Schemes and Validation Study} \label{sec:eval}
\subsection{Numerical Scheme With Boundary Data Extension}\label{sub:qiang_du}
\begin{figure}
    \centering
    \begin{tikzpicture}[scale=3]
\draw[->] (0,0) -- (1.1,0) node[right] {$s$};
\draw[->] (0,0) -- (0,1.1) node[above] {$k(s)$};
\node at (0.7,1.1) {$k(s) \approx \frac{1}{0.0891}{e^{-\frac{1}{(1-x)^2}}}$)};
\draw[thick, blue, domain=0:0.95, samples=100]
    plot (\x,{exp(-1/(1-\x)^2)});

\draw[dashed] (1,0) -- (1,1);
\node[below] at (1,0) {$d=1$};
\end{tikzpicture}

    \caption{Everywhere smooth exponential kernel}
    \label{fig:qiang_du_kernel}
\end{figure}
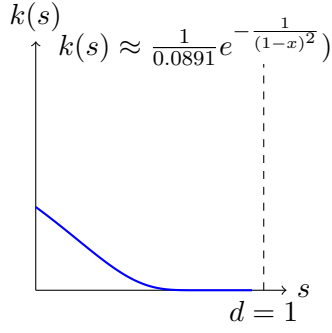

Here we present an adapted Lax--Friedrichs numerical scheme with continuous boundary data extension as in Figure~\ref{fig:Du_bdd}. 
Following the approach in \cite{kachroo2023nonlocal}, we employ an adapted Lax–Friedrichs scheme. The spatial domain 
$(a,b)$ is partitioned into $n$ uniform cells with mesh size  $\Delta x=(b-a)/n$ and time increment $\Delta t$. Let us also define $n_t = T/\Delta t$. 
Then based on \eqref{eq:nonlocal_space_time_conservation_flux}, the time update difference equation is given by,
\begin{align}\label{eq:nlIBVPlwr:NA:LF}
&\rho[n+1,j] \\ \nonumber
&=\frac{1}{2}\big(\rho[n,j+1]+ \rho[n, j-1]) \\ \nonumber
&\quad - \frac{\Delta t}{2 \Delta x}\left[
f(n,j+1,\rho_d[n,j+1])-f(n,j-1,\rho[n,j-1])
\right],
\end{align}
where $\rho_d(t,x)$ is the nonlocal density at the point $(t,x)$ computed via the prescribed look‑ahead convolution.
$\eta(s)$ is the kernel defined in the interval $[0,d]$ where $d<(b-a)$. Since $v_f = 60$, we choose that $\gamma = 0.01$ which is almost equal to $\frac{1}{2 v_f}$. That gives us $\gamma d = 0.01*d$. Next we define $n_d = \frac{d}{\Delta x}$ and $nt_s = \lfloor \frac{\gamma*\Delta x}{\Delta t} \rfloor$. We compute the nonlocal density, as given in \eqref{eq:nonlocal_density_time_space}, using the following equation, 
\begin{equation}
    \begin{gathered}
        \rho_d[n,j] = \bigg(\sum^{n_d-1}_{i=0} \rho[n-i*nt_s,j+i] \eta[i] \bigg) \Delta x,
    \end{gathered}
\end{equation}
Then we use the following kernel expression which satisfies all the conditions given in \cite{du2022space}. We will denote this kernel as smooth exponential kernel (see Figure \ref{fig:qiang_du_kernel}).
\begin{equation}\label{eq:kernel 1}
 \eta(s) = K e^{-\frac{1}{(s-d)^2}}.   
\end{equation}

It can be observed that $\eta(s)$ is $C^1$ and in fact $C^{\infty}$. If we take the derivative we have,
$$ \eta'(s) = \frac{2K  e^{-\frac{1}{(s-d)^2}}}{(s-d)^3}$$ and also satisfies 
$$\eta'(s) \leq - \beta \eta(s)$$
everywhere for $s>0$. Using this kernel and the boundary data, we reconstructed the traffic data for the region $(0, T) \times (a,b)$. If $\rho^{recon}$ is the solution and $\rho^{known}$ is the original data, then we define the error to be
\begin{equation}\label{eq:error_cont}
    er = \frac{\int \int_{(0,T)\times (a,b)} (\rho^{recon}- \rho^{known})^2 dt dx}{\int \int_{(0,T)\times (a,b)} (\rho^{known})^2 dt dx}
\end{equation}
Numerically we compute the error in the following way:
\begin{equation}\label{eq:error_discrete}
    er= \frac{\sum^{n_t}_{i=1}\sum^{n}_{j=1} (\rho^{recon}[i,j]- \rho^{known}[i,j])^2 }{\sum^{n_t}_{i=1}\sum^{n}_{j=1} (\rho^{known}[i,j])^2}
\end{equation}

In addition to the smooth exponential kernel, we also consider several alternative kernel functions. Below, we present two exponential‑type kernels that satisfy the required smoothness properties only on their support. Precisely, we assume
that the kernel to be $C^1{[0,d)}$, supported on $[0,d)$ and $$\eta'(s) \leq -\beta \eta(s), \quad \forall s \in (0,d),$$
for some $\beta >0$. The first one, $\eta(s)$ will be denoted as exponential kernel in the subsequent discussion. We define 

\[\eta(s)=
\begin{cases}
\dfrac{e^{-s/d}}{d\bigl(1-e^{-1}\bigr)}, & 0\le s<d,\\
0, & s\ge d.
\end{cases}\]

\begin{figure}[htbp]
    \hspace{-5mm}\centerline{\includegraphics[width=0.25\textwidth]{ 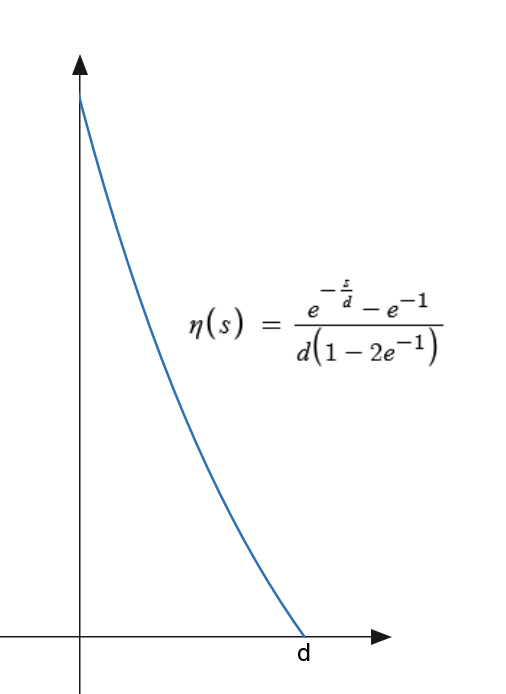}}
    \caption{The Shifted Exponentially Decreasing Kernel
    }
    \label{fig:Shifted exponential_kernel}
\end{figure}

On the other hand, the shifted exponential kernel (shown in Figure \ref{fig:Shifted exponential_kernel}) is defined as 
\[\eta(s)=
\begin{cases}
\dfrac{e^{-s/d}-e^{-1}}{d\bigl(1-2e^{-1}\bigr)}, & 0\le s<d,\\
0, & s\ge d.
\end{cases}\]

Besides, we also run our nonlocal model for the linearly decreasing kernel.

\subsection{Numerical Scheme With Known `Thick' Data}\label{sub:thick_bdd}
Here we present a numerical scheme with known `thick' data as in Figure \ref{fig:thick_bdd}.
 Since part of the solution is already prescribed, the reconstruction error should be evaluated only on the region where the solution is unknown. For this reason, we compute the error as follows:
\begin{equation}\label{eq:error_cont_2}
    er = \frac{\int \int_{(\gamma d,T)\times (a,b-d)} (\rho^{recon}- \rho^{known})^2 dt dx}{\int \int_{(\gamma d,T)\times (a,b-d)} (\rho^{known})^2 dt dx}
\end{equation}
Numerically we compute the error in the following way:
\begin{equation}\label{eq:error_discrete_2}
    er= \frac{\sum^{n_t}_{i=nt_s}\sum^{n-n_d}_{j=1} (\rho^{recon}[i,j]- \rho^{known}[i,j])^2 }{\sum^{n_t}_{i=nt_s}\sum^{n-n_d}_{j=1} (\rho^{known}[i,j])^2}
\end{equation}



\subsection{Case Study and Model Validation}

Next, we present the validation of the proposed nonlocal traffic flow models. We start with the classical case. Then we present the one-dimensional nonlocal (only spatial nonlocal) traffic states, and then transition to two-dimensional (spatial-temporal) settings. The objective is to investigate the influence of (1) kernel length, (2) convolution function, and (3) fixed-length (4) variable length on traffic state estimation accuracy. All experiments are conducted using the trajectory dataset from the Next Generation Simulation (NGSIM) program \cite{ngsim_us101}. The state reconstruction procedure is similar to previous studies \cite {huang2020physics,huang2022physics, huang2023limitations}. The dataset provides detailed vehicle trajectory information on a segment of U.S. Highway 101 in Los Angeles, California. Traffic density fields are reconstructed from the trajectory data using spatial aggregation over fixed grid cells \cite{avila2020data}. The performance of traffic state estimation is evaluated using the relative $L^2$ error between reconstructed density field and the ground truth on US101.

\subsubsection{Classical model}
\begin{figure}[!htbp]
\centering

\subfloat[NGSIM Data\label{fig:ngsim_data}]{
\includegraphics[width=0.48\textwidth]{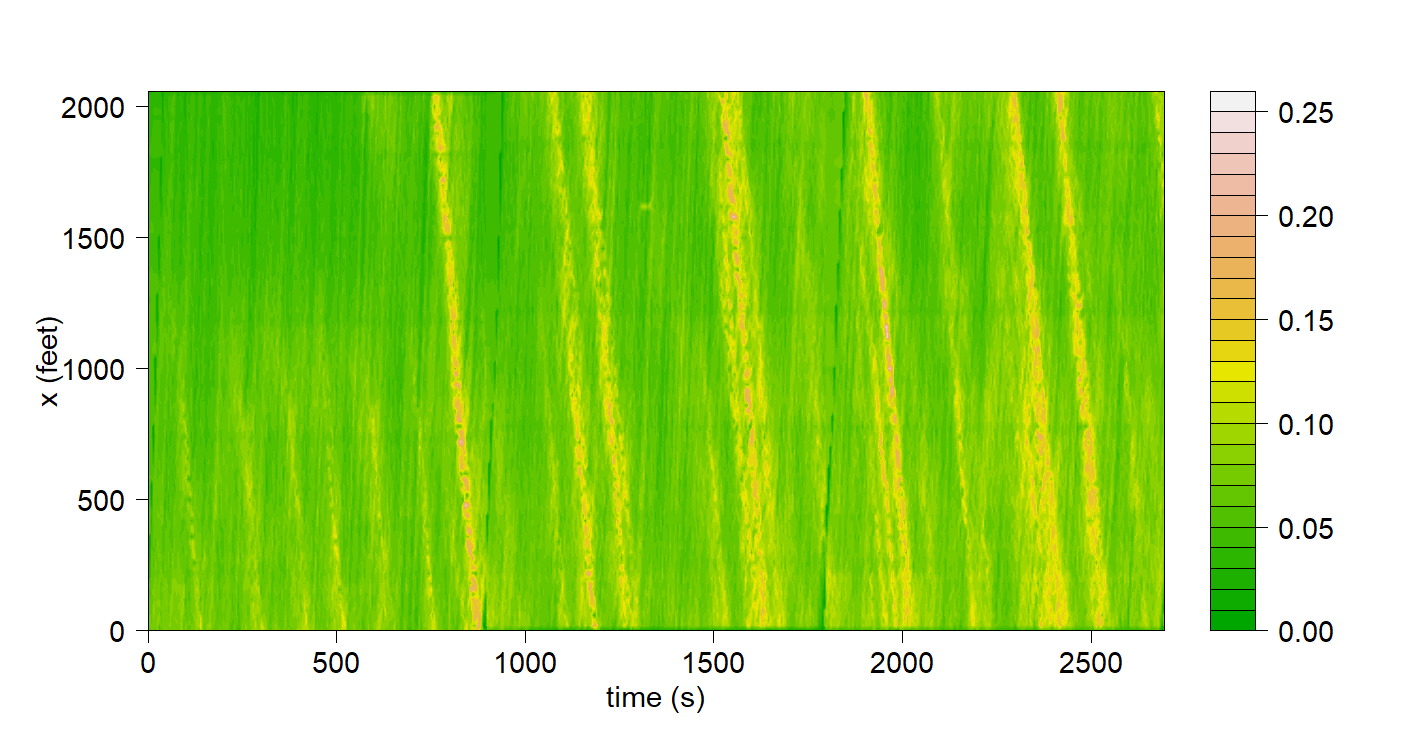}}

\subfloat[Classical Model\label{fig:local}]{
\includegraphics[width=0.48\textwidth]{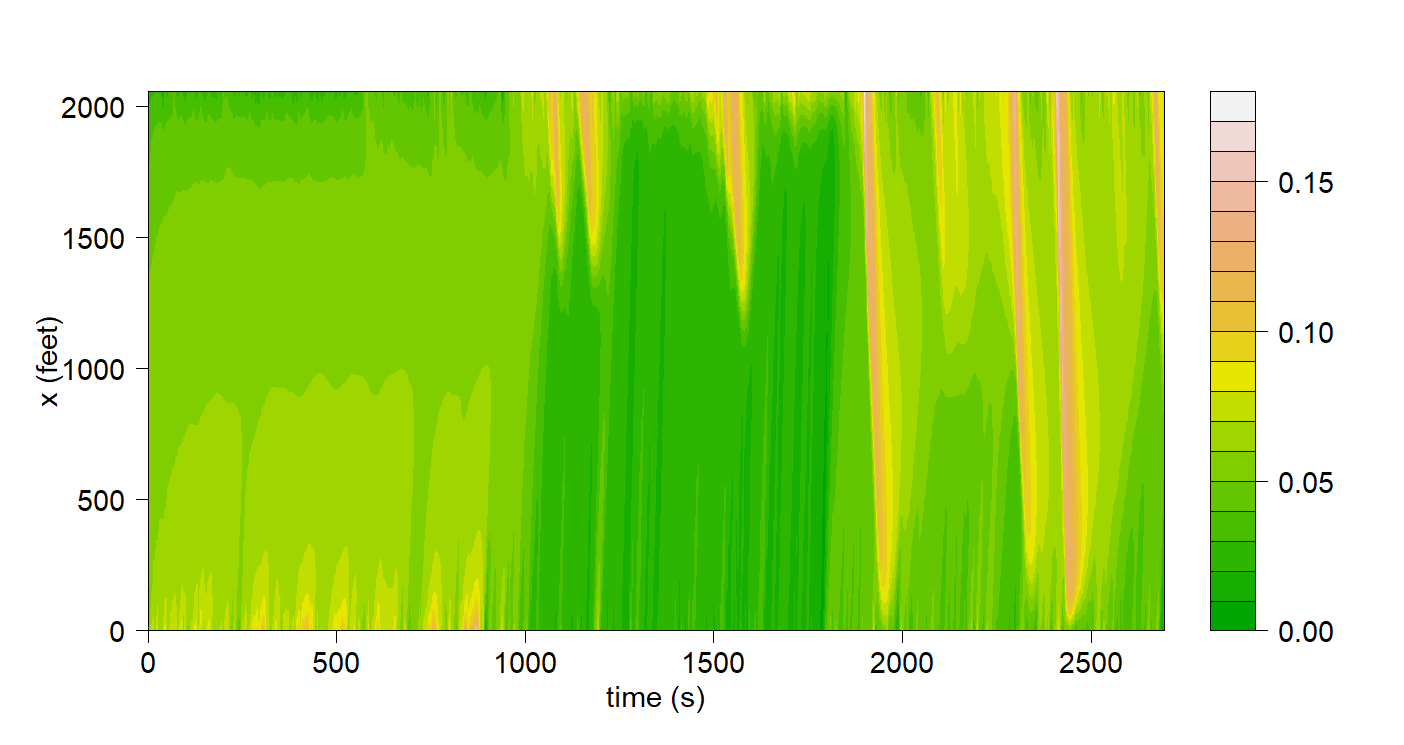}}

\caption{NGSIM Dataset and Vehicle Density Reconstruction with Classical Model}
\label{fig:ngsim}
\end{figure}

We begin with the classical model. In this setting, we consider the boundary value problem in which the initial condition and the boundary value at $x=b$ are prescribed. The corresponding relative error, computed using \eqref{eq:error_discrete}, is 
$0.2133$. The reconstructed solution is shown in Figure\ref{fig:ngsim}(b), while Figure\ref{fig:ngsim}(a) displays the original data.

\subsubsection{1-D Nonlocal Model with Fixed Kernel Size with continuously extended boundary data}

Next we consider the one-dimensional nonlocal vehicle density state with convolution in the spatial domain. For the 1-D reconstructions, we direct the reader's attention to the approach in these papers \cite {9294236, huang2022physics, huang2023limitations}.

\begin{figure*}[htbp]
\centering

\subfloat[Linear kernel\label{fig:0201}]{
\includegraphics[width=0.43\textwidth]{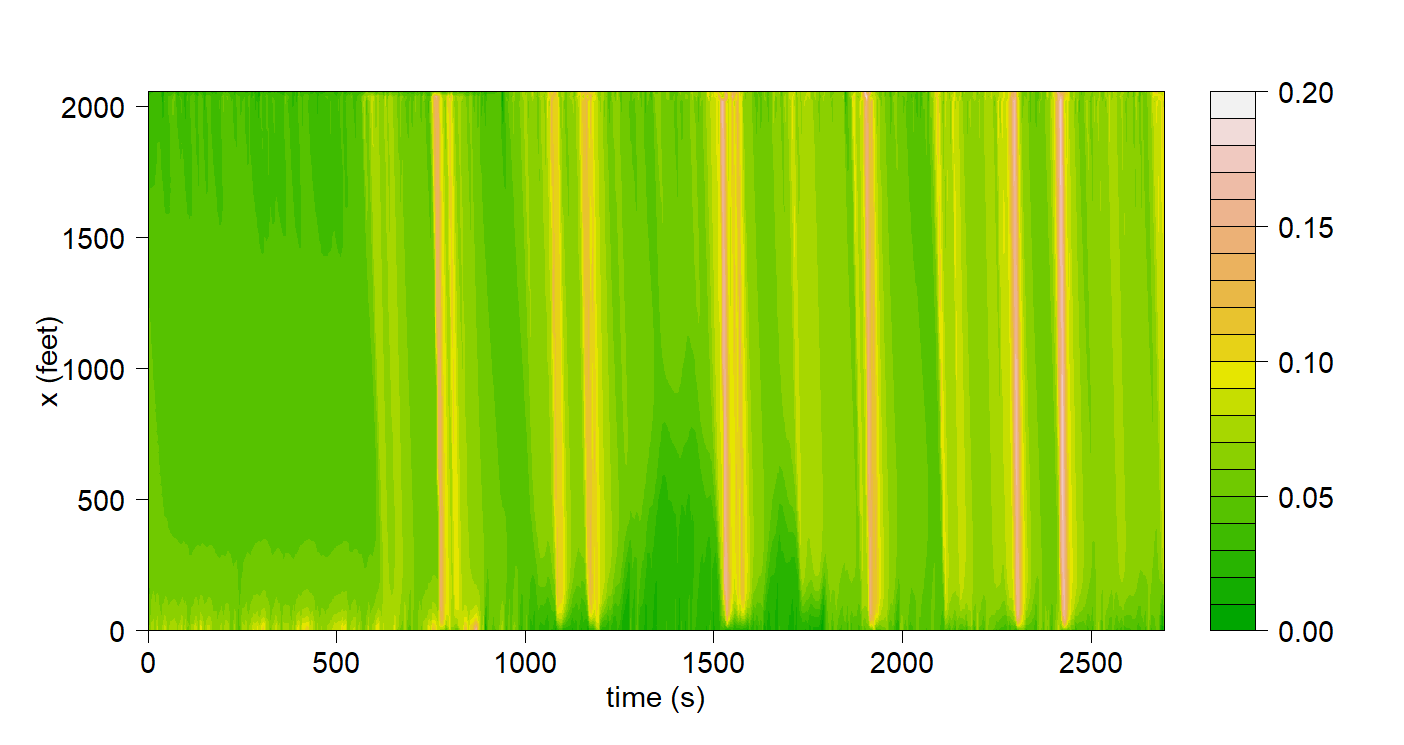}}
\subfloat[Exponential kernel\label{fig:0202}]{
\includegraphics[width=0.43\textwidth]{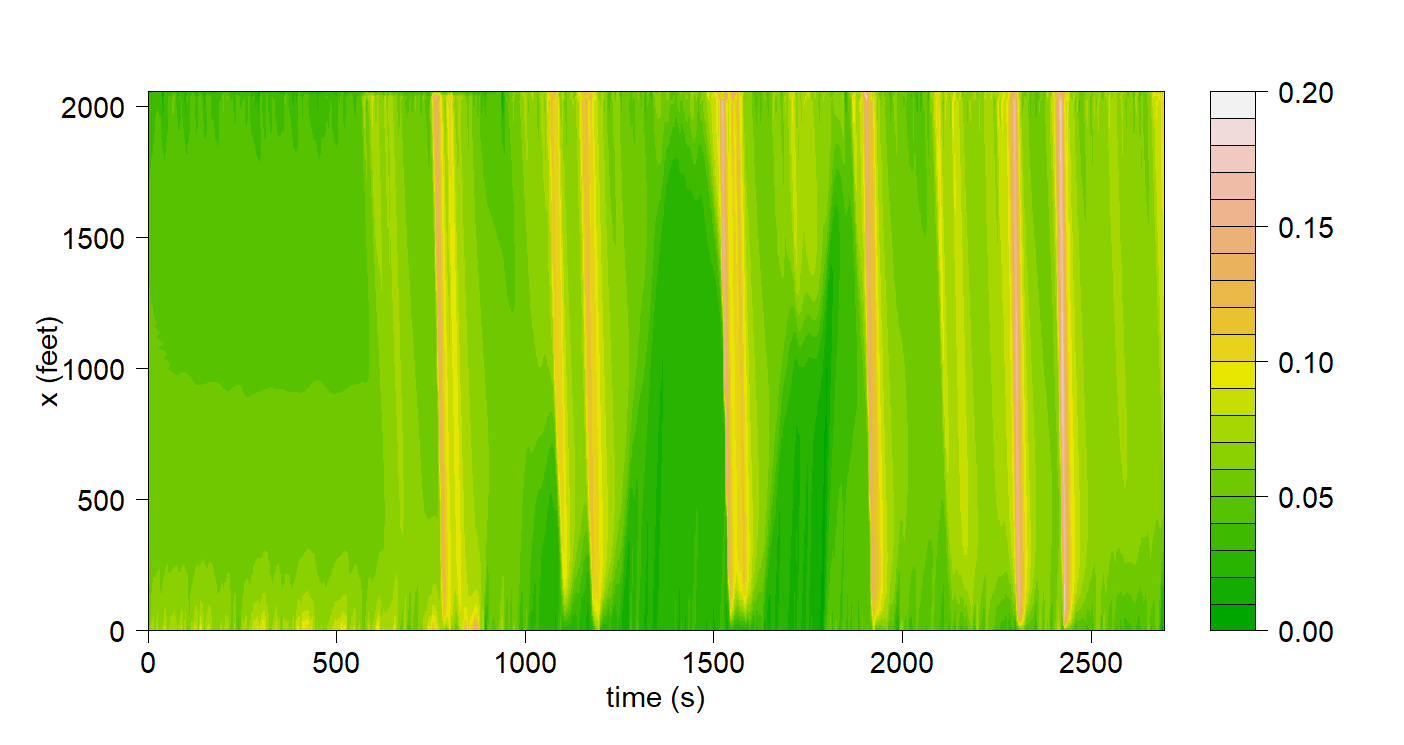}}


\subfloat[Shifted exponential kernel\label{fig:0203}]{
\includegraphics[width=0.43\textwidth]{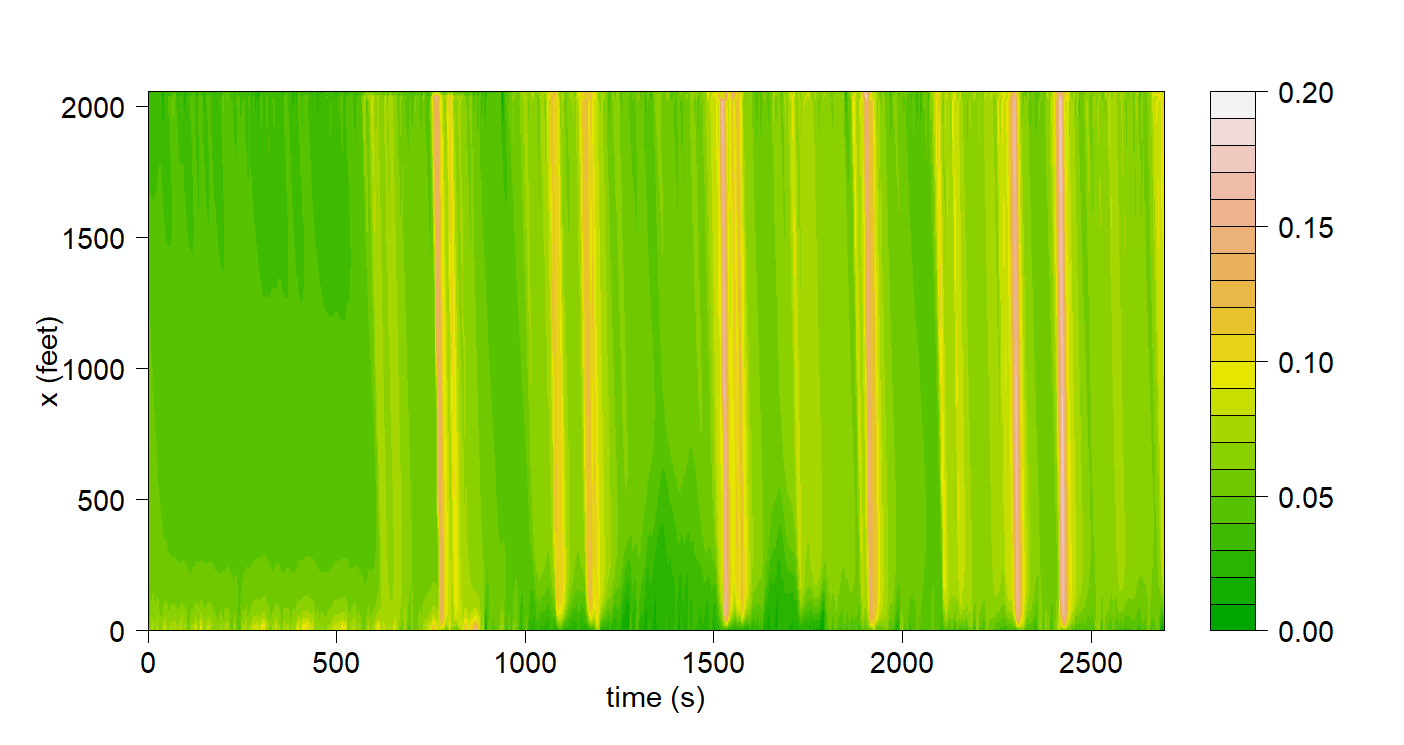}}
\subfloat[Smooth exponential kernel\label{fig:0204}]{
\includegraphics[width=0.43\textwidth]{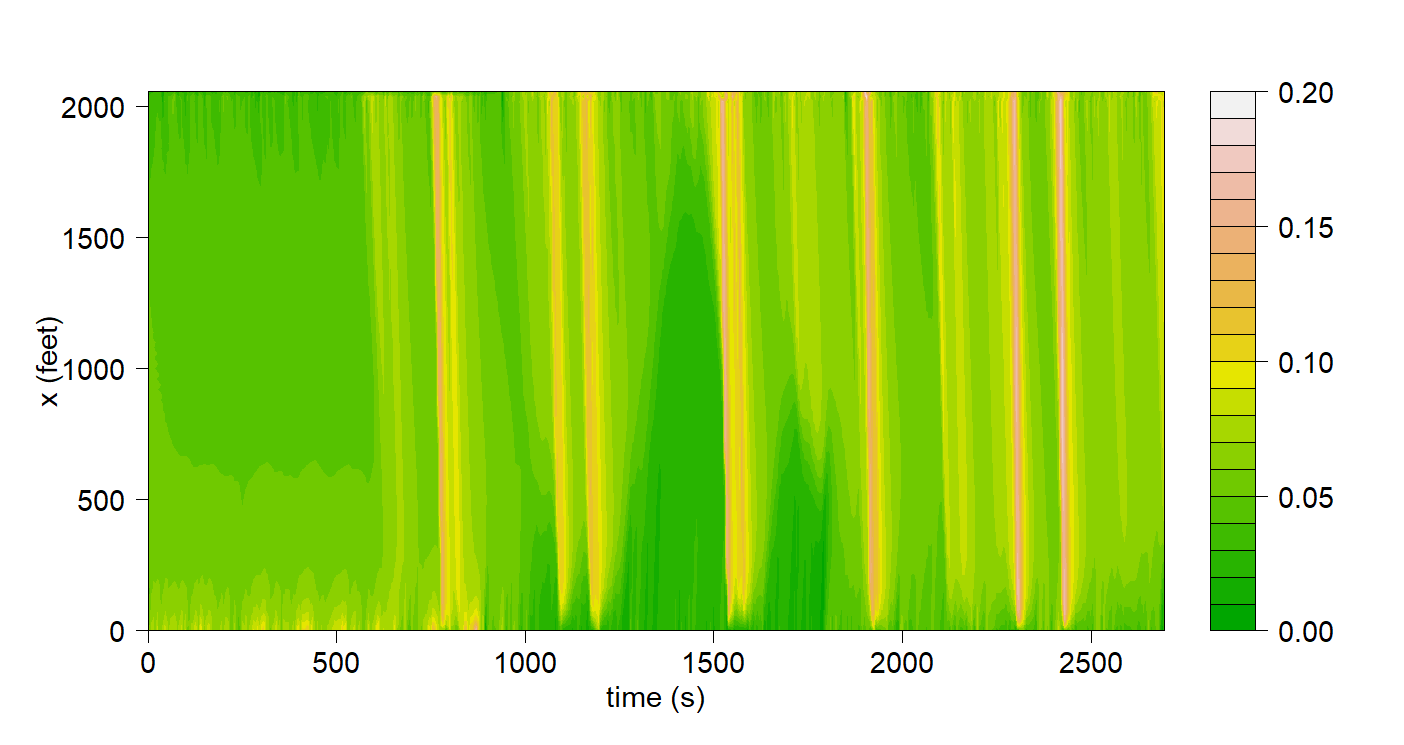}}

\caption{Case I - Vehicle Density Reconstruction with Two-dimensional Nonlocal Models (Kernel Length: 40-ft)}
\label{fig:2d_40ft}

\end{figure*}

\begin{figure*}[h!]
\centering

\subfloat[Linear kernel\label{fig:0301}]{
\includegraphics[width=0.43\textwidth]{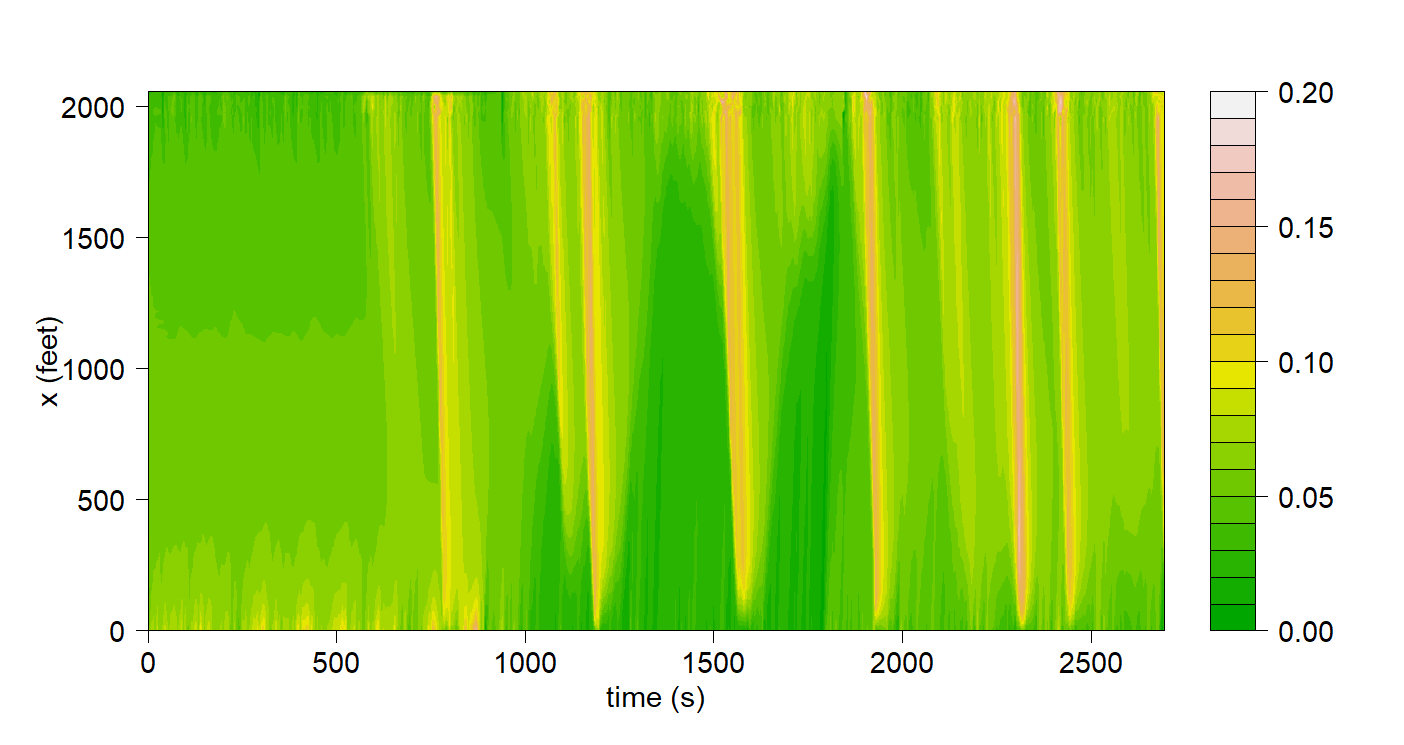}}
\subfloat[Exponential kernel\label{fig:0302}]{
\includegraphics[width=0.43\textwidth]{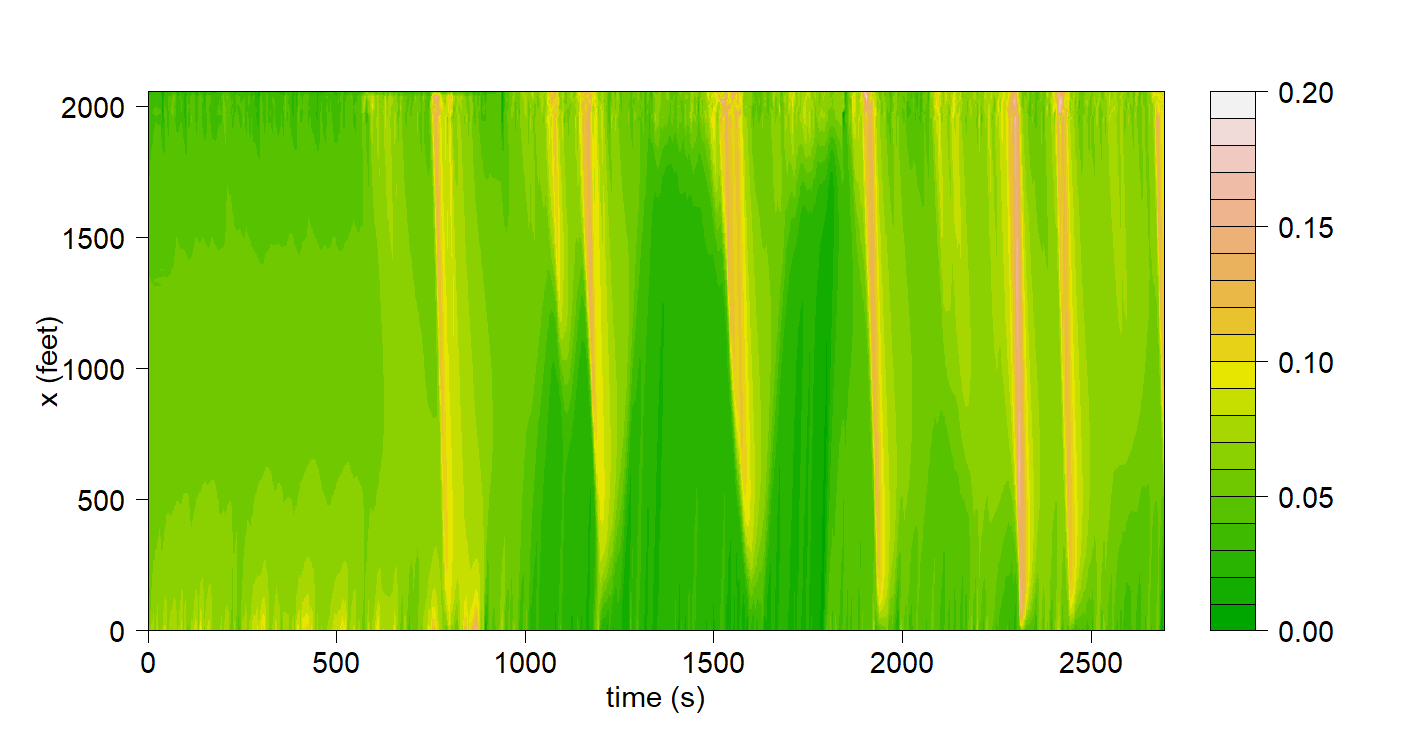}}


\subfloat[Shifted exponential kernel \label{fig:0303}]{
\includegraphics[width=0.43\textwidth]{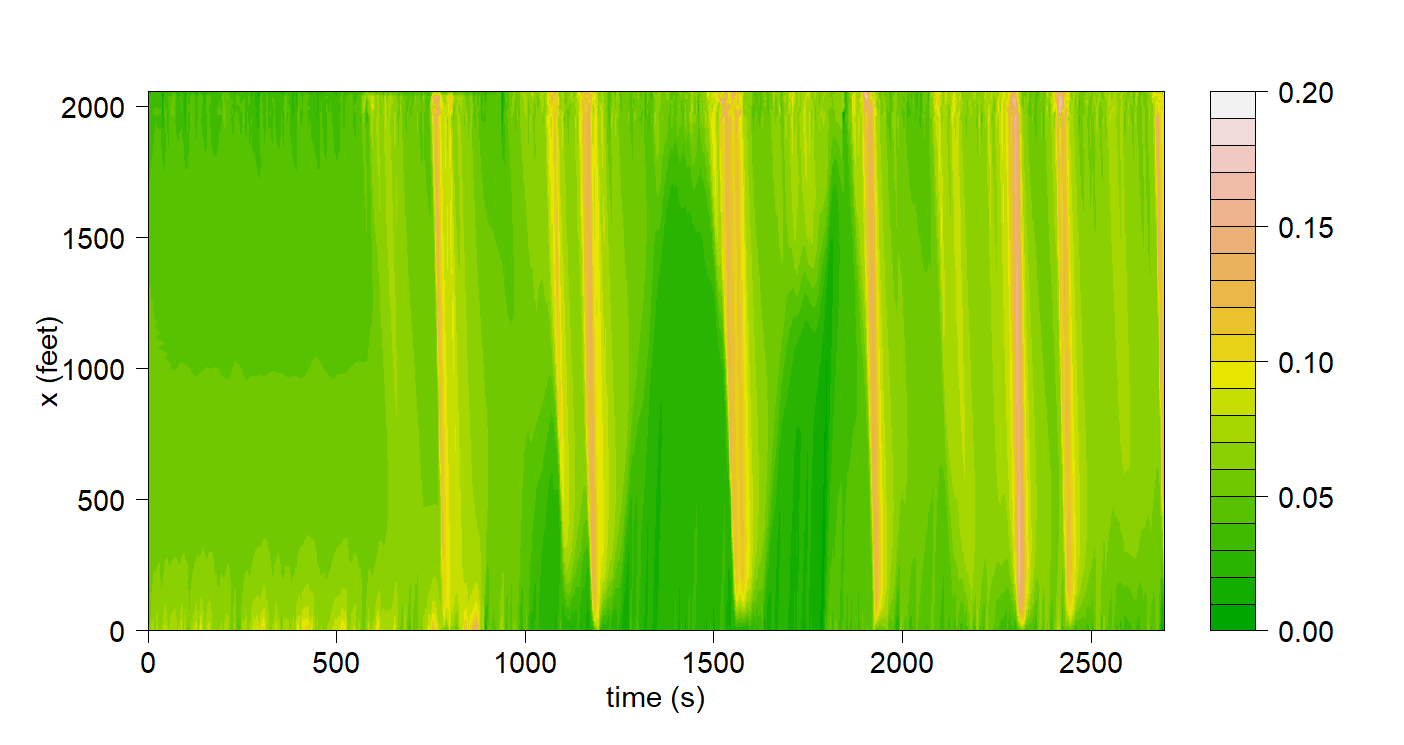}}
\subfloat[Smooth exponential kernel\label{fig:0304}]{
\includegraphics[width=0.43\textwidth]{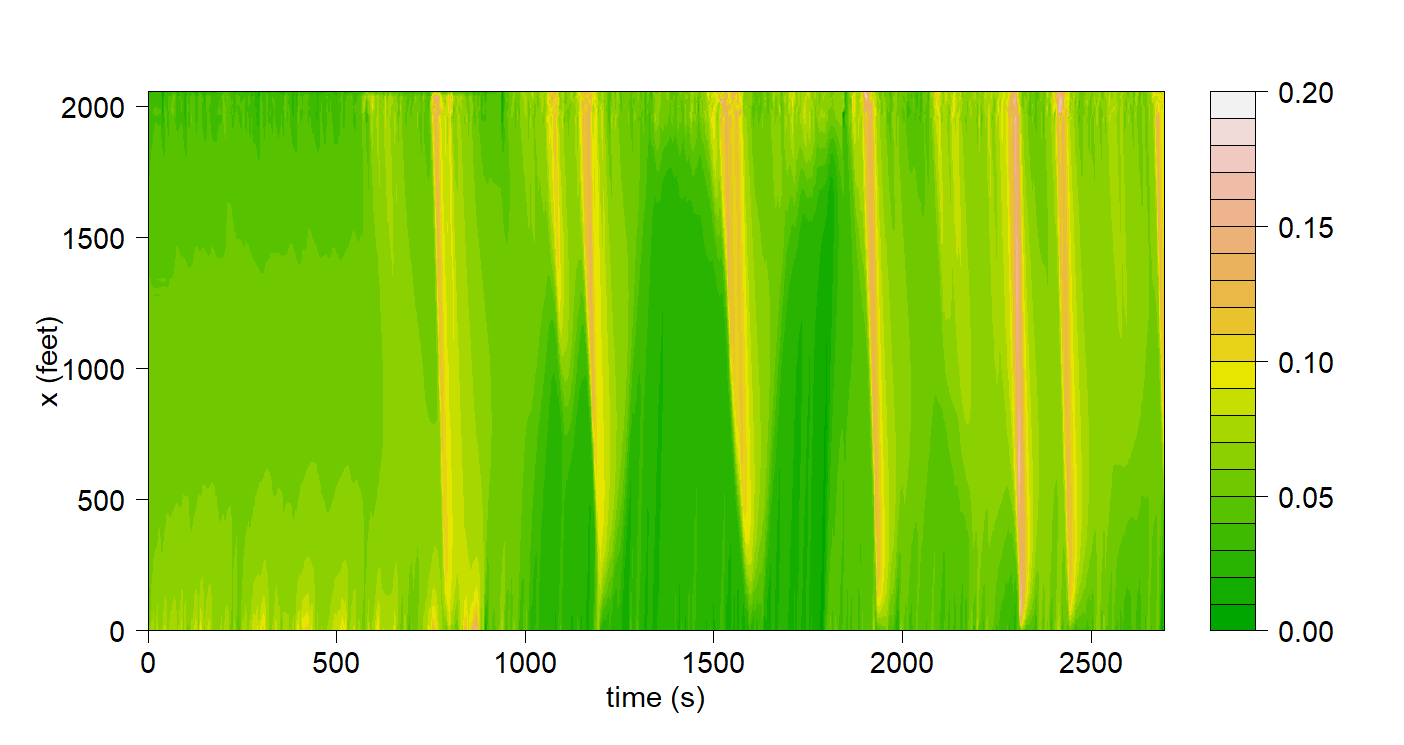}}

\caption{Case II - Vehicle Density Reconstruction with Two-dimensional Nonlocal Models (Kernel Length: 100-ft)}
\label{fig:2d_100ft}

\end{figure*}




\subsubsection{2-D Nonlocal Model with Fixed Kernel Size with continuously extended data}
We now extend the analysis to the two‑dimensional nonlocal traffic model, where the convolution is taken over the space–time domain and the thick boundary data are obtained by continuously extending the thin boundary data, as described in \ref{sub:cont_bdd} and also in \ref{sub:qiang_du}. We investigate that four different kernel configurations:

\begin{itemize}
\item Linear kernel: $K(s) \propto d - s$,
\item Exponential kernel: $K(s) \propto e^{-s/d}$,
\item Shifted exponential kernel: $K(s) \propto e^{-s/d} - e^{-1}$.
\item Smooth exponential kernel: $K(s) \propto \exp\!\left(-\frac{1}{(s-d)^2}\right)$
\end{itemize}

with fixed interaction lengths $$ d = 40 \text{ or } 100 \text{ ft}, $$

All kernels are normalized. We use the $\gamma =0.01$ for our simulation. 
The results are presented in Table~\ref{tab:2d_results_qiang}. We observe that the best case scenario is almost comparable to the corresponding best case in the 1-D case. It is important to note that the boundary and initial value information are same as in the classical case. 

\begin{table}[h]
\centering
\caption{Relative errors with 2-D nonlocal models}
\begin{tabular}{c|c|c|c|c}
Kernel Size (ft) & Linear & Exp. & Shifted Exp. & Smooth Exp. \\
\hline
40 & 0.1424 & 0.1557 & \textbf{0.1423} & 0.1482 \\
100 & 0.1583 & 0.1725 & 0.1564 & 0.1693 \\
\end{tabular}
\label{tab:2d_results_qiang}
\end{table}

\subsubsection{2-D Nonlocal Model with Variable Length Kernels}
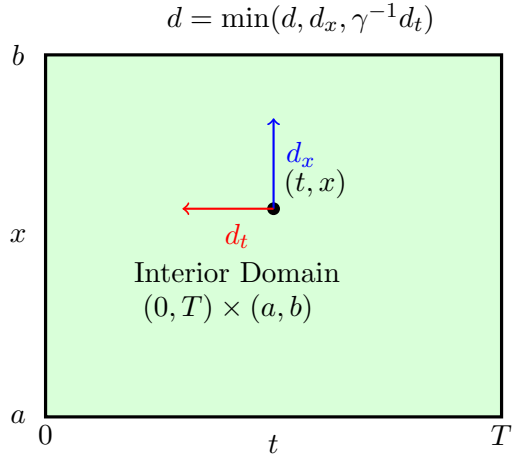
\begin{figure}[h]
\centering
\begin{tikzpicture}[scale=1.2]

\def\T{5}
\def\a{0}
\def\b{4}

\fill[green!15] (0,\a) rectangle (\T,\b);
\draw[very thick] (0,\a) rectangle (\T,\b);


\node at (2.1,1.6) {$\text{Interior Domain}$};
\node at (2.0,1.2) {$(0, T)\times(a,b)$};

\node at (-0.3,2) {$x$};
\node at (2.5,-0.3) {$t$};

\node at (-0.3,0) {$a$};
\node at (-0.3,\b) {$b$};
\node at (0,-0.2) {$0$};
\node at (\T,-0.2) {$T$};

\node at (2.8,4.4) {$d=\min(d, d_x,\gamma^{-1}d_t)$};

\fill (2.5,2.3) circle (2pt);
\node[above right] at (2.5,2.3) {$(t,x)$};

\draw[blue,thick,->] (2.5,2.3) -- (2.5,3.3);
\node[blue] at (2.8,2.9) {$d_x$};

\draw[red,thick,->] (2.5,2.3) -- (1.5,2.3);
\node[red] at (2.1,2.0) {$d_t$};


\end{tikzpicture}
\caption{Variable kernel nonlocal case with local boundary}
\label{fig:thin_bdd}
\end{figure}

Before presenting our best-case results with known `thick' data, we introduce an additional approach for solving the nonlocal boundary value problem using classical boundary and initial conditions. In this setting, the boundary or initial data coincide with those in the classical formulation, or equivalently, with the nonlocal formulation under continuous data extension. This approach is known as the variable-length kernel method, in which the kernel horizon is not fixed a priori but instead varies in both space and time. 
At each spatial location and time step, the admissible kernel length is constrained simultaneously by the pre-determined maximum kernel size, the amount of available past time, and the remaining distance to the right boundary (downstream). 
Precisely, one can use
\begin{equation}
 d(t,x) = \min \{ d, d_x, \gamma^{-1} d_t \}
\end{equation}
where $d$ is the maximum length of admissible kernel, $d_x$ is the distance from the $x=b$ and $d_t$ is the distance from $t=0$, see Figure \ref{fig:thin_bdd}.
As a result, the nonlocal interaction gradually expands as time progresses, while naturally shrinks near the right boundary. The results are shown in Table~\ref{tab:2d_results_variable}.

\begin{table}[h]
\centering
\caption{Relative errors with 2-D nonlocal models, Variable Length Kernel}
\begin{tabular}{c|c|c|c|c}
Max Kernel Size (ft) & Linear & Exp. & Shifted Exp.  & Smooth Exp. \\
\hline
40 & 0.1383  & 0.1539 & \textbf{0.1393} & 0.1446 \\
100 & 0.1544 & 0.1600 &0.1503  & 0.1691 \\
\end{tabular}
\label{tab:2d_results_variable}
\end{table}

Notice that when the maximum allowable kernel length increases to $100$-ft, the deterioration in accuracy observed in the fixed-length case is substantially mitigated. In the variable-length formulation, the effective kernel length is constrained by the past time and the distance to the downstream boundary, preventing excessive averaging when insufficient temporal or spatial information is available. Consequently, the model avoids the over-smoothing effects associated with a uniformly large $100$-ft kernel.

\subsubsection{2-D Nonlocal Model with Fixed Kernel Size with known data}
We now assume that a `thick' initial condition and a `thick' spatial boundary condition are known. This model is described in \ref{sub:known_bdd} and also in \ref{sub:thick_bdd}.
The results are presented in Table~\ref{tab:2d_results}. Again, the shifted exponential kernel 1 with 40-ft length achieves the lowest error, while the standard exponential kernel yields the largest error.

\begin{table}[h]
\centering
\caption{Relative errors with 2-D nonlocal models}
\begin{tabular}{c|c|c|c|c}
Kernel Size (ft) & Linear & Exp. & Shifted Exp.  & Smooth Exp. \\
\hline
40 & 0.1312 & 0.1498 & \textbf{0.1297} & 0.1400 \\
100 & 0.1506 & 0.1645 & 0.1475 & 0.1628 \\
\end{tabular}
\label{tab:2d_results}
\end{table}

The consistently superior performance of the shifted exponential kernel suggests that it offers a balanced representation between local traffic states and the anticipational driving behavior of nonlocal models. Overall, the findings demonstrate that both kernel function and length  influence model accuracy, with shorter-range shifted exponential interactions yielding the most reliable performance.
Figures~\ref{fig:2d_40ft} and \ref{fig:2d_100ft} illustrate the spatial-temporal density reconstruction for all four kernel functions. 

\vspace{-0.15 in}

\subsection{Discussion on results}
The numerical results clearly demonstrate the advantages of nonlocal formulations over the classical LWR model. Among the three nonlocal approaches considered, the 2‑D nonlocal model with known thick data achieves the highest accuracy, with a relative error of 0.1297. The 2‑D nonlocal model with variable‑length kernels performs nearly as well, attaining a relative error of 0.1393, followed closely by the 2‑D nonlocal model with a fixed kernel and continuously extended boundary data, which yields a relative error of 0.1423. In contrast, the classical LWR model exhibits a substantially larger error of 0.2133. These comparisons highlight the clear benefits of using nonlocal traffic flow model, even when only local boundary information is available.

\section{Conclusion} \label{sec:conc}
This work develops and analyzes a new framework for spatial‑temporal nonlocal traffic dynamics, integrating both spatial and temporal interactions into a unified macroscopic flow model. We showed several key analytical properties, including the qualitative behavior under different nonlocal boundary conditions. On the computational side, we designed a numerical scheme tailored to the time‑dependent nonlocal structure and carried out a detailed empirical validation using NGSIM trajectory data. These experiments consistently showed that spatial‑temporal nonlocal models offer substantial improvements over local formulations. The comparison of fixed‑kernel and adaptive‑kernel strategies further clarifies how kernel design influences model performance, offering practical guidance for real‑world implementation. For example, when reliable `thick' boundary data are available, employing a fixed‑kernel nonlocal approach is both natural and effective. However, when such enriched boundary information is absent, one must either revert to the classical model or adopt a nonlocal strategy that compensates for the missing data. One option is to maintain a fixed kernel length and continuously extend the boundary data. The alternative is to use a variable‑length kernel whose support adapts to the distance of each point from the temporal and spatial boundaries. 

While the results demonstrate the promise of nonlocal formulations, several open questions remain. The selection and learning of optimal kernel structures, the treatment of noisy or sparse boundary data, and the extension of the theory to multi‑lane or heterogeneous traffic remain important directions for future research. More broadly, the framework developed here provides a foundation for integrating nonlocal PDE theory with data‑driven methods, opening the door to hybrid modeling approaches that combine mathematical structure with modern machine‑learning tools.

\bibliographystyle{IEEEtran}
\bibliography{references}


\end{document}